\documentclass[12pt]{amsart}
\usepackage{mathrsfs}

\usepackage{amsmath}
\usepackage{amssymb}
\usepackage{amsbsy}
\usepackage{graphicx}
\usepackage{amsfonts}

\newtheorem{lem}{Lemma}[section]%
\newtheorem{theorem}[lem]{Theorem}%
\newtheorem{exam}[lem]{Example}%
\newtheorem{prop}[lem]{Proposition}%

\def\a{\alpha}

 \def\lg{\langle} \def\rg{\rangle}
\def\nd{\mathrel{\bigm|\kern-.7em/}}

\def\f{\noindent}

\def\B{\mathcal{B}}

\def\Aut{\hbox{\rm Aut}}

\def\Cay{\hbox{\rm Cay}}

\def\mod{\hbox{\rm mod }}

\def\demo{\f {\bf Proof.}\hskip10pt}

\def\mz{{\mathbb Z}}

\begin{document}

\title{A family of tetravalent one-regular graphs}
\author[M.~Ghasemi]{Mohsen Ghasemi}
\address{Mohsen Ghasemi,
Department of Mathematics, Urmia University,\newline Urmia  57135,
Iran}\email{m.ghasemi@urmia.ac.ir}

\author[R.~Varmazyar]{Rezvan Varmazyar}
\address{Rezvan Varmazyar,
Department of Mathematics, Azad University, Khoy Branch, Khoy, Iran
}\email{varmazyar@iaukhoy.ac.ir}


\begin{abstract}
A graph is one-regular if its automorphism group acts regularly on
the set of its arcs. In this paper, $4$-valent one-regular graphs of
order $5p^2$, where $p$ is a prime, are classified.
\end{abstract}
\maketitle
\section{Introduction}

In this paper we consider undirected finite connected graphs without
loops or multiple edges. For a graph $X$ we use $V(X)$, $E(X)$,
$A(X)$ and $\Aut(X)$ to denote its vertex set, edge set, arc set and
its full automorphism group, respectively. For $u, v\in V(X)$, $\{u,
v\}$ is the edge incident to $u$ and $v$ in $X$, and $N(u)$ is the
neighborhood of $u$ in $X$, that is, the set of vertices adjacent to
$u$ in $X$. A graph $X$ is said to be {\em vertex-transitive} and
{\em arc-transitive} (or {\em symmetric}) if $\Aut(X)$ acts
transitively on $V(X)$ and $A(X)$, respectively. In particular, if
$\Aut(X)$ acts regularly on $A(X)$, then $X$ is said to be {\em
one-regular}.

Clearly, a one-regular graph  is connected, and it is of valency $2$
 if and only if it is a cycle. In this sense the first non-trivial case is that of cubic
graphs. The first example of a cubic one-regular graph was
constructed by Frucht ~\cite{F} and later on lot of works have been
done along this line (as part of the more general investigation of
cubic arc-transitive graphs) see~\cite{FK2, FK3, FK4, FKW}.
$4$-valent one-regular graphs have also received considerable
attention. In~\cite{C}, $4$-valent one-regular graphs of prime order
were constructed. In ~\cite{M.}, an infinite family of $4$-valent
one-regular Cayley graphs on alternating groups is given. $4$-valent
one-regular circulant graphs were classified in ~\cite{Xu} and
$4$-valent one-regular Cayley graphs on abelian groups were
classified in ~\cite{Xu.j}. Next, one may deduce a classification of
$4$-valent one-regular Cayley graphs on dihedral groups from
~\cite{K.O,WXu,WZh}.   Let $p$ and $q$ be primes. Then, clearly
every $4$-valent one-regular graph of order $p$ is a circulant
graph. Also, by ~\cite{CO, P.W, P.X, WX1, Xu, Xu.j} every $4$-valent
one-regular graph of order $pq$ or $p^2$ is a circulant graph.
Furthermore, the classification of $4$-valent one-regular graphs of
order $3p^2$, $4p^2$, $6p^2$ and $2pq$ are given in ~\cite{FKDZ,Gh1,
GS, JXZYQF}. Along this line the aim of this paper is to classify
$4$-valent one-regular graphs of order $5p^2$, see
Theorem~\ref{main}.

\section{Preliminaries}

In this section, we introduce some notations and definitions as well
as some preliminary results which will be used later in the paper.

For a regular graph $X$, use $d(X)$ to represent the valency of $X$,
and for any subset $B$ of $V(X)$, the subgraph of $X$ induced by $B$
will be denoted by $X[B]$. Let $X$ be a connected vertex-transitive
graph, and let $G\leq \Aut(X)$ be vertex-transitive on $X$. For a
$G$-invariant partition $\B$ of $V(X)$, the {\em quotient graph}
$X_\B$ is defined as the graph with vertex set $\B$ such that, for
any two vertices $B,C\in \B$, $B$ is adjacent to $C$ if and only if
there exist $u\in B$ and $v\in C$ which are adjacent in $X$. Let $N$
be a normal subgroup of $G$. Then the set $\B$ of orbits of $N$ in
$V(X)$ is a $G$-invariant partition of $V(X)$. In this case, the
symbol $X_\B$ will be replaced by $X_N$.

For a positive integer $n$, denote by $\mz_n$ the cyclic group of
order $n$ as well as the ring of integers modulo $n$, by $\mz_n^*$
the multiplicative group of $\mz_n$ consisting of numbers coprime to
$n$, by $D_{2n}$ the dihedral group of order $2n$, and by $C_n$ and
$K_n$ the cycle and the complete graph of order $n$, respectively.
We call $C_n$ an {\em $n$-cycle}.

For  a finite group $G$ and a subset $S$ of $G$ such that $1 \notin
S$ and $S=S^{-1}$, the {\em Cayley graph} $\Cay(G,S)$ on $G$ with
respect to $S$ is defined to have vertex set $G$ and edge set
$\{\{g,sg\} \mid g \in G, s \in S\}$. Given a $g \in G$, define the
permutation $R(g)$ on $G$ by $x\mapsto xg$, $x \in G$. The
permutation group $R(G)=\{R(g) \mid g \in G\}$ on $G$ is called the
{\em right regular representation} of $G$. It is easy to see that
$R(G)$ is isomorphic to $G$, and it is a regular subgroup of the
automorphism group $\Aut(\Cay(G,S))$. Also it is easy to see that
$X$ is connected if and only if $G=\langle S \rangle$, that is, $S$
is a connection set. Furthermore, the group $\Aut(G,S)=\{\alpha \in
\Aut(G) \mid S^\alpha=S\}$ is a subgroup of $\Aut(\Cay(G,S))$.
Actually, $\Aut(G,S)$ is a subgroup of $\Aut(\Cay(G,S))_1$, the
stabilizer of the vertex $1$ in $\Aut(\Cay(G,S))$. A Cayley graph
$\Cay(G,S)$ is said to be {\em normal} if $R(G)$ is normal in
$\Aut(\Cay(G,S))$. Xu \cite{Xu2}, proved that $\Cay(G,S)$ is normal
if and only if $\Aut(\Cay(G,S))_1=\Aut(G,S)$. Suppose that $\alpha
\in \Aut(G)$. One may easily prove that $\Cay(G,S)$ is normal if and
only if $\Cay(G,S^\alpha)$ is normal. Also later much subsequent work
was done along this line (see \cite {BFSX, FX, GhZ, WXu,  WZh}). 

For $u\in V(X)$, denote by $N_X(u)$ the {\it neighbourhood} of $u$
in $X$, that is, the set of vertices adjacent to $u$ in $X$. A graph
$\widetilde{X}$ is called a {\it covering} of a graph $X$ with
projection $p:\widetilde{X}\rightarrow X$ if there is a surjection
$p:V{(\widetilde{X})}\rightarrow V(X)$ such that
$p|_{N_{\widetilde{X}}({\tilde{v}})}:{N_{\widetilde{X}}({\tilde{v}})}\rightarrow
{N_{X}(v)}$ is a bijection for any vertex $v\in V(X)$ and
$\tilde{v}\in p^{-1}(v)$. A covering $\widetilde{X}$ of $X$ with a
projection $p$ is said to be {\it regular} (or {\it K-covering}) if
there is a semiregular subgroup $K$ of the automorphism group
Aut($\widetilde{X}$) such that graph $X$ is isomorphic to the
quotient graph $\widetilde{X}/K$, say by $h$, and the quotient map
$\widetilde{X}\rightarrow \widetilde{X}/K$ is the composition $ph$
of $p$ and $h$ (for the purpose of this paper, all functions are
composed from left to right). If $K$ is cyclic or elementary abelian
then $\widetilde{X}$ is called a {\it cyclic} or an {\it elementary
abelian covering} of $X$, and if $\widetilde{X}$ is connected $K$
becomes the covering transformation group. In this case we also say
$p$ is a {\it regular covering projection}.  The {\it fibre} of an
edge or a vertex is its preimage under $p$. An automorphism of
$\widetilde{X}$ is said to be {\it fibre-preserving} if it maps a
fibre to a fibre, while every covering transformation maps a fibre
on to itself. All of fibre-preserving automorphisms form a group
called the {\it fibre-preserving group}.

Let $\widetilde{X}$ be a $K$-covering of $X$ with a projection $p$.
If $\alpha$$\in$ Aut($X$) and $\widetilde{\alpha}$$\in$
Aut($\widetilde{X}$) satisfy $\widetilde{\alpha}p=p\alpha$, we call
$\widetilde{\alpha}$ a {\it lift} of $\alpha$, and $\alpha$ the {\it
projection} of $\widetilde{\alpha}$. Concepts such as a lift of a
subgroup of Aut($X$)  and the projection of a subgroup of
Aut$(\widetilde{X})$ are self-explanatory. The lifts and the
projections of such subgroups are of course subgroups in
Aut($\widetilde{X}$) and Aut($X$) respectively.

For two groups $M$ and $N$, $N\rtimes M$ denotes a semidirect
product of $N$ by $M$. For a subgroup $H$ of a group $G$, denote by
$C_G(H)$ the centralizer of $H$ in $G$ and by $N_G(H)$ the
normalizer of $H$ in $G$. Then $C_G(H)$ is normal in $N_G(H)$.

\begin{prop} {\rm\cite[Chapter I, Theorem~4.5]{Hup}}\ \  \label{NC2}
The quotient group $N_G(H)/C_G(H)$ is isomorphic to a subgroup of
the automorphism group $\Aut(H)$ of $H$.
\end{prop}

Let $G$ be a permutation group on a set $\Omega$ and $\a\in \Omega$.
Denote by $G_\a$ the stabilizer of $\a$ in $G$, that is, the
subgroup of $G$ fixing the point $\a$. We say that $G$ is {\em
semiregular} on $\Omega$ if $G_\a=1$ for every $\a\in \Omega$ and
{\em regular} if $G$ is transitive and semiregular. For any $g\in
G$, $g$ is said to be {\em semiregular} if $\lg g\rg$ is
semiregular.

\begin{prop} {\rm\cite[Chapter I, Theorem~4.5]{HW}}\ \  \label{NC1}
Every transitive abelian group $G$ on a set $\Omega$ is regular.
\end{prop}



The following proposition is due to Praeger et al, refer to
\cite[Theorem~1.1]{GP}.

\begin{prop}\label{reduction}
Let $X$ be a connected $4$-valent $(G,1)$-arc-transitive graph. For
each normal subgroup $N$ of $G$, one of the following holds:
\begin{enumerate}
\item [{\rm (1)}] $N$ is transitive on $V(X)$;

\item [{\rm (2)}] $X$ is bipartite and $N$ acts transitively on
each part of the bipartition;

\item [{\rm (3)}] $N$ has $r\geq 3$ orbits on $V(X)$, the quotient
graph $X_N$ is a cycle of length $r$, and $G$ induces the full
automorphism group $D_{2r}$ on $X_N$;

\item [{\rm (4)}] $N$ has $r\geq 5$ orbits on $V(X)$, $N$ acts
semiregularly on $V(X)$, the quotient graph $X_N$ is a connected
$4$-valent $G/N$-symmetric graph, and $X$ is a $G$-normal cover of
$X_N$.
\end{enumerate}
Moreover, if $X$ is also $(G,2)$-arc-transitive, then case $(3)$ can
not happen.
\end{prop}



The following classical result is due to Wielandt \cite[
Theorem~3.4]{HW}

\begin{prop}  \label{NC3}
Let $p$ be a prime and let $P$ be a Sylow $p$-subgroup of a
permutation group $G$ acting on a set $\Omega$. Let $ \omega \in
\Omega$. If $p^m$ divides the length of the $G$-orbit containing
$\omega$, then $p^m$ also divides the length of the $P$-orbit
containing $\omega$.
\end{prop}

To state the next result we need to introduce a family of $4$-valent
graphs that were first defined in ~\cite{GP1}. The graph
$C^{\pm1}(p;5p,1)$ is  defined to have the vertex set
$\mathbb{Z}_p\times \mathbb{Z}_{5p}$ and edge set
$\{(i,j)(i\pm1,j+1)| i \in\mathbb{Z}_{p}, j \in \mathbb{Z}_{5p}\}$.
Also  from~\cite[Definition~$2.2$]{GP1}, the graphs $C^{\pm
1}(p;5p,1)$ are Cayley graphs over $\mathbb{Z}_p\times
\mathbb{Z}_{5p}$ with connection set
$\{(1,1),(-1,1),(-1,-1),(1,-1)\}.$ In the proof of
Theorem~\ref{main}, we will need $C^{\pm 1}(p;5p,1)$ with $p> 11$.
It can be readily checked from~\cite[Definition~$2.2$]{GP1} that for
$p> 11$ these graphs are actually normal Cayley graphs over
$\mathbb{Z}_p\times \mathbb{Z}_{5p}$.

\begin{prop}{\rm\cite[Theorem 1.1]{GP1}} \label{gp-sym}
Let $X$ be a connected, $G$-symmetric, $4$-valent graph of order
$5p^2$, and let $N=\mathbb{Z}_p$ be a minimal normal subgroup of $G$
with orbits of size $p$, where $p$ is an odd prime. Let $K$ denote
the kernel of the action of $G$ on $V(X_N)$. If $X_N=C_{5p}$ and
$K_v\cong \mathbb{Z}_2$ then $X$ is  isomorphic to
$C^{\pm1}(p;5p,1)$.
\end{prop}

The graphs defined in  ~\cite[Lemma 8.4]{GP1} are all one-regular
(see ~\cite[Section~$8$]{GP1}) and therefore we refer to ~\cite{GP1}

for an intrinsic description of these families.

\begin{prop}{\rm\cite[Theorem 1.2]{GP1}}\label{gp-sym1}
Let $X$ be a connected, $G$-symmetric, $4$-valent graph of order
$5p^2$, and let $N=\mathbb{Z}_p \times \mathbb{Z}_p$ be a minimal
normal subgroup of $G$ with orbits of size $p^2$, where $p$ is an
odd prime. Let $K$ denote the kernel of the action of $G$ on
$V(X_N)$. If $X_N=C_{5}$ and $K_v\cong \mathbb{Z}_2$ then $X$ is
isomorphic to  one of the graphs in ~\cite[Lemma 8.4]{GP1}.
\end{prop}

Finally in the following example we introduce $G(5p;2,2,u)$, which
first was defined in ~\cite{P.W}.

\begin{exam}\label{examp1}
Let $2$ be a divisor of $p-1$. Let $H(5,2)=\langle a \rangle$, let
$t \in \mathbb{Z}_{p}^{*}$ be such that $t \in -H(p,2)$, and let $u$
be the least common multiple of $2$ and the order of $t$ in
$\mathbb{Z}_{p}^{*}$. Then $X=G(5p; 2,2, u)$ is defined as the graph
with vertex set $$V(X)=\mathbb{Z}_5\times \mathbb{Z}_p=\{(i,x)| i
\in \mathbb{Z}_5, x \in \mathbb{Z}_p\}$$ such that vertices $(i,x)$
and $(j,y)$ are adjacent if and only if there is an integer $l$ such
that $j-i=a^l$ and $y-x \in t^lH(p,2)$. Also $X$ as defined above is
independent of the choice of generator $a$ of $H(5,2)$ up to
isomorphism, and $X$ is also independent of the choice of $t$, such
that $lcm\{o(t),2\}=u$, up to isomorphism. Moreover, the above graph
is circulant, that is, admits a cyclic group of automorphisms of
order $5p$ acting regularly on vertices.

\end{exam}

We may extract the following results from ~\cite[pp. 76-80]{CF}.

\begin{prop}\label{class}
Let $p$ be a prime and $p>5$. Also let $G$ be a non-abelian group of
order $5p^2$.
\begin{enumerate}
\item [{\rm (i)}] If $G$ has a normal subgroup of order $p$, say $N$, such that $G/N$
is cyclic, then $G$ is isomorphic to $\langle x,y, z|
x^p=y^5=z^p=[x,z]=[y,z]=1, y^{-1}xy=x^i \rangle$, where $i^5\equiv
1$ $(\mod p)$ and  $(i,p)=1$;

\item [{\rm (ii)}] If $G$ has a normal subgroup of order $p^2$, say $N$, such that $G/N$
is cyclic, then $G$ is isomorphic to $\langle x,y| x^{p^2}=y^5=1,
y^{-1}xy=x^i \rangle$, where $i^5\equiv 1$ $(\mod p^2)$.

\end{enumerate}

\end{prop}

\section{one-regular graphs of order $5p^2$}

For proving the main theorem we need  the following two lemmas.

\begin{lem}\label{lemma-a}
Let $p$ ba a prime, $p> 5$ and $G=\langle x,y, z|
x^p=y^5=z^p=[x,z]=[y,z]=1, y^{-1}xy=x^i\rangle$, where $i^5\equiv 1$
$(\mod p)$ and $(i,p)=1$. Then there is no $4$-valent one-regular
normal Cayley graph $X$ of order $5p^2$ on $G$.
\end{lem}
\demo Suppose to the contrary that $X$ is a $4$-valent one-regular
normal Cayley graph Cay$(G,S)$ on $G$ with respect to the generating
set $S$. Since $X$ is one-regular and normal, the stabilizer
$A_1=\Aut(G,S)$ of the vertex $1 \in G$ is transitive on $S$ and so
that elements in $S$ are all of the same order. The elements of $G$
of order $5$ lie in $\langle x, y \rangle$ and the elements of $G$
of order $p$ lie in $\langle x,z \rangle$. Since $X$ is connected,
$G=\langle S \rangle$ and hence $S$ consists of elements of order
$5p$. Denote by $\mathcal{S}_{5p}$ the elements of $G$ of order
$5p$. Therefore
$$S\subseteq \mathcal{S}_{5p}=\{x^sy^tz^j | s \in
\mathbb{Z}_{p}, t \in \mathbb{Z}_5^{*}, j \in \mathbb{Z}_p^{*}\}.$$
 Clearly $\sigma:
x\mapsto x^s, y\mapsto y, z\mapsto z^j$ $(s, j \neq 0)$ is an automorphism of $G$,
we may suppose that $$S=\{xy^tz, y^{-t}x^{-1}z^{-1}, x^my^nz^k,
y^{-n}x^{-m}z^{-k}\}.$$  Since $\Aut(G,S)$ acts transitively on $S$,
it implies that there is $\alpha \in \Aut(G,S)$ such that
$(xy^tz)^\alpha=y^{-t}x^{-1}z^{-1}$. Since $[x,z]=1$, and $[y,z]=1$,
the element $z^\alpha$ needs to commute with $x^\alpha$ and
$y^\alpha$. Thus $(xz)^\alpha
(y^t)^\alpha=y^{-t}x^{-1}z^{-1}=x^{-i^{-4t}}y^{-t}z^{-1}$. We may
assume that $(y^t)^\alpha=x^{t_1}y^{t_2}$, where $t_1 \in
\mathbb{Z}_p$, and $t_2 \in \mathbb{Z}_5^{*}$. Thus $(xz)^\alpha
x^{t_1}y^{t_2}=x^{-i^{-4t}}y^{-t}z^{-1}$ and so
$(xz)^\alpha=x^{-i^{-4t}}x^{-t_1i^{4(-t-t_2)}}y^{-t-t_2}z^{-1}$.
Since $o(xz)=p$, we have $t=-t_2$. Therefore
$(xz)^\alpha=x^{-i^{-4t}-t_1}z^{-1}$. Also let
$z^\alpha=x^{s_1}z^{s_2}$ where $s_1, s_2 \in \mathbb{Z}_p$. So
$(x)^\alpha=x^{-i^{-4t}-s_1-t_1}z^{-1-s_2}$. Since $z^\alpha$
commutes with $(y^t)^\alpha$, it follows that $s_1=0$ or
$i^{4t_2}=1$.

Since  $t_2 \in \mathbb{Z}_5^{*}$  and  $i^5\equiv 1$ $(\mod p)$, it
follows that $i^{4t_2}\neq 1$. Thus we may  suppose that $s_1=0$.
Therefore $x^\alpha=x^{-i^{-4t}-t_1}z^{-1-s_2}$,
$(y^t)^\alpha=x^{t_1}y^{-t}$, $z^\alpha=z^{s_2}$. Since
$x^{y^t}=x^{i^t}$, we have
$(x^\alpha)^{(y^t)^\alpha}=(x^\alpha)^{i^t}$ and so $s_2=-1$ and
$(-i^{-4t}-t_1)(i^{-4t_2}-i^t)=0$.  Since $t=-t_2$ and $t_2 \in
\mathbb{Z}_5^{*}$, we have $(i^{-4t_2}-i^t)\neq 0$. Thus we may
suppose that  $(-i^{-4t}-t_1)=0$. Therefore
$x^\alpha=z^{-1-s_2}=z^{0}=1$, a contradiction.\\

\begin{lem}\label{lemma-b}
Let $p$ ba a prime, $p> 5$ and $G=\langle x,y| x^{p^2}=y^5=1,
y^{-1}xy=x^i\rangle$, where $i^5\equiv 1$ $(\mod p^2)$. Then there
is no $4$-valent one-regular normal Cayley graph $X$ of order $5p^2$
on $G$.
\end{lem}
\demo Suppose to the contrary that $X$ is a $4$-valent one-regular
normal Cayley graph Cay$(G,S)$ on $G$ with respect to the generating
set $S$. Since $X$ is one-regular and normal, the stabilizer
$A_1=\Aut(G,S)$ of the vertex $1 \in G$ is transitive on $S$ and so
that elements in $S$ are all of the same order. Clearly $x^p$ is the
only element of order $p$. Also $x^r$ where $r \in
\mathbb{Z}_{p^2}^{*}$ are the only elements of order $p^2$.  The
elements of $G$ of order $5$ lie in $\langle x, y \rangle$. Since
$X$ is connected, $G=\langle S \rangle$ and hence $S$ consists of
elements of order $5$. Denote by $\mathcal{S}_{5}$ the elements of
$G$ of order $5$. Therefore
$$S\subseteq \mathcal{S}_{5}=\{x^ry^s | r \in \mathbb{Z}_{p^2}, s \in \mathbb{Z}_{5}^{*} \}.$$
 Clearly $\sigma:
x\mapsto x^r, y\mapsto y$ $(r \neq 0)$ is an automorphism of $G$, we may suppose
that $S=\{xy^s, y^{-s}x^{-1}, x^uy^v, y^{-v}x^{-u}\}.$ Since
$\Aut(G,S)$ acts transitively on $S$, it implies that there is
$\alpha \in \Aut(G,S)$ such that $(xy^s)^\alpha=y^{-s}x^{-1}$. We
may assume that $y^\alpha=x^{m}y^{n}$, where $m \in
\mathbb{Z}_{p^2}$, $n \in \mathbb{Z}_{5}^{*}$. Also let
$x^\alpha=x^r$, where $r  \in \mathbb{Z}_{p^2}^{*}$. Therefore
$x^r(x^my^n)^s=y^{-s}x^{-1}$, and so $ns=-s$. Thus $s=0$ or $n=-1$.
Clearly, $s\neq 0$, and so $n=-1$. Now $y^\alpha=x^{m}y^{-1}$. Since
$x^y=x^{i}$, we have $(x^\alpha)^{y^\alpha}=(x^\alpha)^{i}$ and so
$ri^4-ri=0$. Thus $i^3=1$, a contradiction.\\

 Let $X$ be a tetravalent one-regular graph of order $5p^2$.
 If $p\leq 11$, then $|V(X)|=20$, $45$, $125$, $245$, or
 $605$. Now, a complete census of the tetravalent arc-transitive graphs
 of order at most $640$ has been recently obtained by Poto\v cnik, Spiga and Verret~\cite{spiga,spiga1}.
 Therefore, a quick inspection through this
 list (with the invaluable help of \texttt{magma} (see ~\cite{BosCanPla})) gives the number of tetravalent
 one-regular graphs in the case that $p\leq 11$. Thus we may suppose that $p>11$.

The following result is the main result of this paper.

\begin{theorem}\label{main}
Let $p$ be a prime. A $4$-valent graph $X$ of order $5p^2$ is
$1$-regular if and only if one of the following holds:
\begin{description}
\item[(i)] $X$ is a Cayley graph over $\langle x,y| x^p=y^{5p}=[x,y]=1\rangle$,  with connection sets $\{y,y^{-1},
 xy, x^{-1}y^{-1}\}$ and $\{y,y^{-2}, xy, x^{-2}y^{-2}\}$;
\item[(ii)] $X$ is  connected arc-transitive circulant graph with respect to every connection set $S$;
\item[(iii)] $X$ is one of the graphs described in~\cite[Lemma 8.4
]{GP1}.
\end{description}
\end{theorem}

\demo Let $X$ be a $4$-valent one-regular graph of order $5p^2$.
 If $p\leq 11$, then $|V(X)|=20, 45$, $125$, $245$, or  $605$.
Now, a complete census of the $4$-valent arc-transitive graphs
 of order at most $640$ has been recently obtained by Poto\v cnik, Spiga and Verret~\cite{spiga,spiga1}.
 Therefore, a quick inspection through this list (with the invaluable help of \texttt{magma}) gives the proof of the theorem in the case that $p\leq 11$.

Now, suppose that $p> 11$. Let $A=\Aut(X)$ and let $A_v$ be the
stabilizer in $A$ of the vertex $v\in V(X)$. Let $P$ be a Sylow
$p$-subgroup of $A$. Since $A$ is one-regular, it follows that
$|A|=20p^2$.  We show that $P$ is normal in $A$. Since $|A|=20p^2$,
the Sylow's theorems show that the number of Sylow $p$-subgroups of
$A$ is equal to $|A:N_A(P)|=1+kp$, for some $k\geq 0$. If $k=0$,
then $P$ is normal in $A$ and thus we may assume that $k\geq 1$.
Now, $1+kp$ divides $20$ and this is possible if and only if $k=1$
and $p=19$. Now $|A:N_A(P)|=20$. So $N_A(P)=P$ and $C_A(P)=N_A(P)$.
  Therefore, by the Burnside's $p$-complement theorem~\cite[page~$76$]{Wer},
  we see that $A$ has a normal subgroup $N$ of order $20$. In particular, $P$ acts by conjugation as a group
  of automorphisms on $N$.  As a group of order $20$ does not
  admit non-trivial automorphisms of order $19$, we see that $P$ centralizes $A$.
  Thus $A\cong N\times P$ and $P$ is normal in $A$.\\

Assume first that $P$ is cyclic. Let $X_P$ be the quotient graph
of $X$ relative to the orbits of $P$ and  let $K$ be the kernel
of $A$ acting on $V(X_P)$. By Proposition~\ref{NC3}, the orbits
of $P$ are of length $p^2$. Thus $|V(X_P)|=5$, $P\leq K$ and
$A/K$ acts arc-transitively on $X_P$. By
Proposition~\ref{reduction}, either $X_P\cong C_5$ and hence
$A/K\cong D_{10}$ forcing that $|K|=2p^2$, or $P$ acts
semiregularly on V$(X)$, the quotient graph $X_P$ is a
tetravalent connected $A/P$-arc-transitive graph and $X$ is a
regular cover of $X_P$. First assume that  $X_P\cong C_{5}$. If
$A/P$ is an abelian then, since $A/K$ is  a quotient group of
$A/P$, also $A/K$ is an abelian. But since $A/K$ is
vertex-transitive on $X_P$, Proposition~\ref{NC1}, implies that
it is regular on $X_P$, contradicting arc-transitivity of $A/K$
on $X_P$. Thus $A/P$ is non-abelian group. Clearly $K$ is not
semiregular on V$(X)$. Then $K_v \cong \mathbb{Z}_2$, where $v
\in V(X)$. By Proposition~\ref{NC2}, $A/C \lesssim
\mathbb{Z}_{p(p-1)}$, where $C=C_A(P)$. Since $A/P$ is not
abelian we have that $P$ is a proper subgroup of $C$. If $C\cap
K\neq P$, then $C\cap K=K$ ($|K|=2p^2$). Since $K_v$ is a Sylow
$2$-subgroup of $K$, $K_v$ is characteristic in $K$ and so normal
in $A$, implying that $K_v=1$, a contradiction. Thus $C\cap K=P$
and $1\neq C/P=C/C\cap K \cong CK/K\unlhd A/K\cong D_{10}$. If
$C/P\cong \mathbb{Z}_2$, then $C/P$ is in the center of $A/P$ and
since $A/P/C/P\cong A/C$ is cyclic, $A/P$ is abelian, a
contradiction. It follows that $|C/P| \in \{5,10\}$, and hence
$C/P$ has a characteristic subgroup of order $5$, say $H/P$. Thus
$|H|=5p^2$ and $H/P\unlhd A/P$, implies that $H\unlhd A$. In
addition since $H\leq C=C_A(P)$, we have that $H$ is abelian.
Clearly $|H_v| \in \{1,5\}$. If $|H_v|=5$, then $H_v$ is a Sylow
$5$-subgroup of $H$, implying that $H_v$ is characteristic in $H$.
The normality of $H$ in $A$ implies that $H_v \unlhd A$, forcing
$H_v=1$, a contradiction. If $H_v=1$, then since $|H|=5p^2$, $H$
is regular on V$(X)$. It follows that $X$ is a Cayley graph on an
abelian group with a cyclic Sylow $p$-subgroup $P$. By elementary
group theory, we know that up to isomorphism $\mathbb{Z}_{5p^2}$,
where $p>11$, is the only abelian group with a cyclic Sylow
$p$-subgroup. Also by ~\cite[Theorem 7]{Xu}, $X$ is one-regular.

Now assume  that $X_P$ is a tetravalent connected $A/P$-symmetric
graph. Clearly, $X_P\cong K_5$ and by ~\cite[Theorem 2.2]{BFSX},
$X_P$ is non-normal Cayley graph on $\mathbb{Z}_5$. On the other
hand $A/P$ is isomorphic to a subgroup of index $6$ in
$\Aut(K_5)\cong S_5$. Thus $A/P$ is isomorphic to affine group
AGL$(1,5)=\mathbb{Z}_5 \rtimes \mathbb{Z}_4$. Therefore $A/P$ has
a normal subgroup of order $5$, say $PM/P$. Thus $PM\unlhd A$ and
$PM$ is transitive on V$(X)$. Since $|PM|=5p^2$, $PM$ is also
regular on V$(X)$, implying that $X$ is a normal Cayley graph on
$PM$. If $PM$ is an abelian group, then  $PM$ is isomorphic to
$\mathbb{Z}_{5p^2}$. Also if $PM$ is not abelian, then by
Proposition~\ref{class} part (ii), $PM$ is isomorphic to $\langle
x,y| x^{p^2}=y^5=1, y^{-1}xy=x^i \rangle$, where $i^5\equiv 1$
$(\mod p^2)$. If $PM\cong \mathbb{Z}_{5p^2}$, then by
~\cite[Theorem 7]{Xu} $X$ is one-regular.   In the latter case,
$X$ is not one-regular, by Lemma~\ref{lemma-b}.


Now assume that $P$ is an elementary abelian. Suppose first that
$P$ is a minimal normal subgroup of $A$, and consider the
quotient graph $X_P$ of $X$ relative to the orbits of $P$. Let
$K$ be the kernel of $A$ acting on V$(X_P)$. By
Proposition~\ref{reduction}, either $X_P\cong C_5$ and hence
$A/K\cong D_{10}$ forcing that $|K|=2p^2$, or $P$ acts
semiregularly on V$(X)$, the quotient graph $X_P$ is a
tetravalent connected $A/K$-arc-transitive graph and $X$ is a
regular cover of $X_P$. First assume that  $X_P\cong C_{5}$. Thus
$K_v=\mathbb{Z}_2$. Proposition~\ref{gp-sym1}  implies that $X$ is
isomorphic to one of the graphs described
in~\cite[Lemma~$8.4$]{GP1}.

Now assume that $X_P$ is a tetravalent connected $A/P$-symmetric
graph. So $X$ is a $\mathbb{Z}_p\times \mathbb{Z}_p$-regular
 cover of $K_5$. By ~\cite[Table 1]{Kuz}, AGL$(1,5)$, lifts along
 $p$. Now we use the fact that the lift of an $s$-regular group that
 lifts along a regular covering projection is $s$-regular (see
 ~\cite{Djo,Djo1}). We recall that AGL$(1,5)$ is a one-regular
 subgroup of $\Aut(K_5)$. Now by ~\cite[Theorem 2.1, Propositions 3.4,
 3.5]{Kuz}, $X$ is not one-regular.

Suppose now that $P$ is not a minimal normal subgroup of $A$.
Then a minimal normal subgroup $N$ of $A$ is isomorphic to
$\mathbb{Z}_p$. Let $X_N$ be  the quotient graph of $X$ relative
to the orbits of $N$ and let $K$ be the kernel of $A$ acting on
V$(X_N)$. Then $N\leq K$ and $A/K$ is transitive on V$(X_N)$,
moreover we have that $|V(X_N)|=5p$. By
Proposition~\ref{reduction}, $X_N$ is a cycle of length $5p$, or
$N$ acts semiregularly on V$(X)$, the quotient graph $X_N$ is
$4$-valent connected $A/N$-arc-transitive graph and $X$ is a
regular cover of $X_N$. If $X_N\cong C_{5p}$, and hence $A/K\cong
D_{10p}$, then $|K|=2p$ and thus $K_v\cong \mathbb{Z}_2$. Applying
Proposition~\ref{gp-sym}, we get that $X$ is  isomorphic to
$C^{\pm1}(p;5p,1)$. If, however $X_N$ is a $4$-valent connected
$A/N$-symmetric graph, then, by Proposition~\ref{reduction}, $X$
is a covering graph of a symmetric graph of order $5p$. By
~\cite{P.W}, $G(5p;2,2,u)$ is the just $4$-valent symmetric graph
of order $5p$ (see Example~\ref{examp1}). Observe that in this
case a one-regular subgroup of automorphism contains a normal
regular subgroup isomorphic to $\mathbb{Z}_5 \times
\mathbb{Z}_p$. Let $H$ be a one-regular subgroup of automorphism
of $X_N$. Since $X$ is one-regular graph, $A$ is the lift of $H$.
Since $H$ contains a normal regular subgroup isomorphic to
$\mathbb{Z}_5 \times \mathbb{Z}_p$ also $A$ contains a normal
regular subgroup. Therefore $X$ is a normal Cayley graph of order
$5p^2$. Since $A/\mathbb{Z}_p\cong H$ and $\mathbb{Z}_5 \times
\mathbb{Z}_p \unlhd H$, there exists a normal subgroup $G$ of $A$
such that $G/\mathbb{Z}_p\cong \mathbb{Z}_{p}\times
\mathbb{Z}_5$. If $G$ is an abelian group, then  $G$ is
isomorphic to $\mathbb{Z}_p\times \mathbb{Z}_{5p}$, or
$\mathbb{Z}_{5p^2}$. Also if $G$ is not abelian, then  by
Proposition~\ref{class} part (i), $G$ is isomorphic to $\langle
x,y, z| x^p=y^5=z^p=[x,z]=[y,z]=1, y^{-1}xy=x^i \rangle$, where
$i^5\equiv 1$ $(\mod p)$ and $(i,p)=1$. If $G\cong
\mathbb{Z}_{5p^2}$, then by ~\cite[Theorem 7]{Xu} $X$ is
one-regular. Also if  $G\cong \mathbb{Z}_p\times \mathbb{Z}_{5p}$
then by ~\cite[Proposition 3.3, Example 3.2]{Xu.j} $X$ is
isomorphic to either $\Cay(\mathbb{Z}_p\times \mathbb{Z}_{5p},
\{a,a^{-1}, ab, a^{-1}b^{-1}\})$ or $\Cay(\mathbb{Z}_p\times
\mathbb{Z}_{5p}, \{a,a^{-2}, ab, a^{-2}b^{-2}\})$ which are
$1$-regular. These graphs are in
Theorem~\ref{main} part (ii). Finally, in the latter case,  $X$ is not  one-regular,  by Lemma~\ref{lemma-a}. This complete the proof. \\

\end{document}